\documentstyle[12pt]{article}

\input amssym.def
\input amssym.tex

\newcommand{\eop}{\bigstar}  

\newcommand{\card}[1]{{\vert #1 \vert} }

\newcommand{\otp}[1]{\hbox{otp($#1$)}}

\newcommand{\cf}{{\rm cf}}

\newcommand{\implies}{\Longrightarrow}

\newenvironment{Proof}{\noindent{\bf Proof.}}{\par\bigskip} 

\newtheorem{THEOREM}{Theorem}[section]

\newtheorem{Conclusion}[THEOREM]{Conclusion}

\newtheorem{LEMMA}[THEOREM]{Lemma}

\newtheorem{Main Theorem}[THEOREM]{Main Theorem}
\newenvironment{main Theorem}{\begin{Main Theorem}} 
{\end{Main Theorem}}

\newtheorem{Theorem}[THEOREM]{Theorem}

\newtheorem{Definition}[THEOREM]{Definition}

\newtheorem{Conventions}[THEOREM]{Conventions}

\newtheorem{Main Definition}[THEOREM]{Main Definition}
\newenvironment{main definition}{\begin{Main Definition}}
{\end{Main Definition}}

\newtheorem{Lemma}[THEOREM]{Lemma}

\newtheorem{Notation}[THEOREM]{Notation}

\newtheorem{Convention}[THEOREM]{Convention}

\newtheorem{Note}[THEOREM]{Note}

\newtheorem{Observation}[THEOREM]{Observation}

\newtheorem{Remark}[THEOREM]{Remark}

\newtheorem{Main Fact}[THEOREM]{Main Fact}
\newenvironment{main Fact}{\begin{Main Fact}}{\end{Main Fact}}

\newtheorem{Fact}[THEOREM]{Fact}

\newtheorem{Subfact}[THEOREM]{Subfact}

\newtheorem{Claim}[THEOREM]{Claim}

\newtheorem{Main Claim}[THEOREM]{Main Claim}
\newenvironment{main claim}{\begin{Main Claim}}{\end{Main Claim}}

\newtheorem{Corrolary}[THEOREM]{Corrolary}

\newtheorem{Subclaim}[THEOREM]{Subclaim}

\newtheorem{Corollary}[THEOREM]{Corollary}

\newtheorem{Example}[THEOREM]{Example}

\newtheorem{Proposition}[THEOREM]{Proposition}

\newtheorem{Discussion}[THEOREM]{Discussion}

\newenvironment{Proof of the Subfact}
{\noindent{\bf Proof of the Subfact.}}{\par\bigskip}

\newenvironment{Proof of the Theorem}
{\noindent{\bf Proof of the Theorem.}}{\par\bigskip}

\newenvironment{Proof of the Conclusion}
{\noindent{\bf Proof of the Conclusion.}}{\par\bigskip}

\newenvironment{Proof of the Observation}
{\noindent{\bf Proof of the Observation.}}{\par\bigskip}

\newenvironment{Proof of the Fact}
{\noindent{\bf Proof of the Fact.}}{\par\bigskip}

\newenvironment{Proof of the Lemma}
{\noindent{\bf Proof of the Lemma.}}{\par\bigskip}

\newenvironment{Proof of the Claim}
{\noindent{\bf Proof of the Claim.}}{\par\bigskip}

\newenvironment{Proof of the Subclaim}
{\noindent{\bf Proof of the Subclaim.}}{\par\medskip}

\newenvironment{Proof of the Main Claim}
{\noindent{\bf Proof of the Main Claim.}}{\par\bigskip}

\catcode`\@=11
\def\@begintheorem#1#2{\rm \trivlist \item[\hskip \labelsep{\bf #1\ #2.}]}
\def\@opargbegintheorem#1#2#3{\rm \trivlist
      \item[\hskip \labelsep{\bf #1\ #2\ (#3).}]}
\catcode`\@=12

\newcommand{\elementary}{\prec}

\newcommand{\into}{\rightarrow}

\newcommand{\rest}{\upharpoonright}  

\newcommand{\nacc}{\mathop{\rm nacc}}
\newcommand{\acc}{\mathop{\rm acc}}

\newcommand{\gl}{\mathop{\rm gl}}

\newcommand{\deq}{\buildrel{\rm def}\over =}

\newcommand{\HH}{{\cal H}}

\newcommand{\PP}{{\cal P}}

\newcount\skewfactor
\def\mathunderaccent#1#2 {\let\theaccent#1\skewfactor#2
\mathpalette\putaccentunder}
\def\putaccentunder#1#2{\oalign{$#1#2$\crcr\hidewidth
\vbox to.2ex{\hbox{$#1\skew\skewfactor\theaccent{}$}\vss}\hidewidth}}
\def\name{\mathunderaccent\tilde-3 }

\begin{document}

\title{On versions of $\clubsuit$ on cardinals larger than $\aleph_1$}

\author{Mirna D\v zamonja\\
School of Mathematics\\
University of East Anglia\\
Norwich, NR47TJ, UK\\
\scriptsize{M.Dzamonja@uea.ac.uk}\\
Saharon Shelah\\
Mathematics Department\\
Hebrew University of Jerusalem\\
91904 Givat Ram, Israel\\
and\\
Rutgers University\\
New Brunswick, New Jersey\\
USA\\
\scriptsize{shelah@sunset.huji.ac.il}}
\date{}
\maketitle
\begin{abstract}
We give two results on guessing unbounded
subsets of $\lambda^+$. The first is a positive result and applies to the
situation of $\lambda$ regular
and at least equal to $\aleph_3$, while the second is a negative
consistency result which applies to the situation of $\lambda$ a singular
strong limit
with $2^\lambda>\lambda^+$. The first result shows that in $ZFC$
there is a guessing of unbounded subsets of $S^{\lambda^+}_\lambda$.
The second result is a consistency result (assuming a supercompact
cardinal exists) showing that a natural guessing fails. A result of
Shelah in \cite{Sh 667}
shows that if $2^\lambda=\lambda^+$ and $\lambda$ is a strong
limit singular, then the corresponding guessing holds.

Both results are also connected to an earlier result of D\v zamonja-Shelah
in which they showed that a certain version
of $\clubsuit$ holds at a successor of singular just in $ZFC$.
The first result here shows that the result of Fact \ref{545} can to
a certain extent be extended to the successor of a regular. The negative
result here gives limitations to the extent to which one can hope to
extend the mentioned D\v zamonja-Shelah result.

{\footnote{This paper is numbered 685 (10/96) in Saharon Shelah's list of
publications. Both authors thank NSF
for partial support by their grant number NSF-DMS-97-04477, as well
as
the United States-Israel Binational Science Foundation for a partial
support through a BSF grant. Our thanks also go to Ofer Shafir, for
pointing out and correcting a problem in an earlier version.
}}
\end{abstract}

\baselineskip=16pt
\binoppenalty=10000
\relpenalty=10000
\raggedbottom

\section{Introduction and background}
The combinatorial principle $\clubsuit$ is
a weakening of $\diamondsuit$ which, at $\aleph_1$,
means that there is a sequence $\langle
A_\delta:\,\delta\mbox{ limit }<\omega_1\rangle$ such that every $A_\delta$
is an unbounded subset of $\delta$, and for every unbounded subset
$A$ of $\omega_1$, there are stationarily many $\delta$ such that
$A_\delta\subseteq A$. One can weaken this statement in various ways, for
example requiring $\card{A_\delta\setminus A}<\aleph_0$ in place of $A_\delta
\subseteq A$ above, and in general the weakened statements are not
equivalent to $\clubsuit$ (as opposed to the situation with $\diamondsuit$,
see D\v zamonja-Shelah \cite{DjSh 576} and Kunen's \cite{Kunen}
respectively), and are not provable in $ZFC$. The question we 
consider here is if the corresponding situation holds at cardinals
larger than $\aleph_1$. As an example of earlier results
in this direction and connected to the statement of our main
theorem, we mention a
result of Erd\"os, Dushnik and Miller, and a much more recent one of Shelah.
Erd\"os, Dushnik and Miller prove in \cite{DuMi}, that
if $\alpha<\lambda^+$, then $\alpha$ can be written as $\bigcup_{n<\omega}
A_n$, where for each $n$ we have $\otp{A_n}<\lambda^n$
(compare this with Note \ref{Cn} below). 
Shelah proves in \cite{Sh 572} that if $\lambda=\cf(\lambda)>\kappa>\aleph_0$,
then there is a sequence
\[
\langle (C_\delta=\{\alpha_{\delta,\varepsilon}:\,\varepsilon<\lambda\},
h_\delta):\,\delta\in
S^{\lambda^+}_\lambda\rangle
\]
with each $\{\alpha_{\delta,\varepsilon}:\,\varepsilon<\lambda\}$ a continuous
increasing
sequence with $\sup_{\varepsilon<\lambda}\alpha_{\delta,\varepsilon}=\delta$, and
$h_\delta:\,C_\delta\into\kappa$ an onto function, such that for every club $E$
of $\lambda^+$, for stationarily many $\delta<\lambda^+$ we have, for every
$i<\kappa$
\[
\{\varepsilon<\lambda:\,\alpha_{\delta,\varepsilon},\alpha_{\delta,\varepsilon+1}\in
E\,\,\&\,\,h(\alpha_{\delta,\varepsilon})=i\}
\]
is stationary in $\lambda$
(see Notation \ref{oznake} below for $S^{\lambda^+}_\lambda$).
Note that it is
not known if the analogous result
holds when ``$\alpha_{\delta,\varepsilon},
\alpha_{\delta,\varepsilon+1}\in E$" is replaced by
$``\alpha_{\delta,\varepsilon}, \alpha_{\delta,\varepsilon+1},
 \alpha_{\delta,\varepsilon+2}\in E"$.
If it does, it would have interesting consequences regarding the
generalized
Suslin hypothesis, see \cite{KjSh 449}.

We prove that if $\lambda$ is regular and at least equal to $\aleph_3$,
then just in $ZFC$ a version of $\clubsuit$ by which unbounded
subsets of $S^{\lambda^+}_\lambda$ (i.e. the ordinals $<\lambda^+$
of cofinality $\lambda$) are guessed, holds.
If $\lambda\ge \aleph_2, 2^{\aleph_0}$ we obtain a similar version of guessing.
See Theorem \ref{general}.

Another result along these lines is one of D\v zamonja and Shelah
from \cite{DjSh 545}:

\begin{Definition}\label{guesso} Suppose that $\lambda$ is a cardinal.
$\clubsuit^\ast_{-\lambda}(\lambda^+)$ is the statement saying that there is
a sequence $\langle \PP_\delta:\,\delta\mbox{ limit }<\lambda^+\rangle$
such that
\begin{description}
\item{(i)} $\PP_\delta$ is a family of $\le\card{\delta}$ unbounded subsets
of $\delta$,
\item{(ii)} For $a\in \PP_\delta$ we have $\otp{a}<\lambda$,
\item{(iii)} For all $X\in [\lambda^+]^{\lambda^+}$, there is a club $C$ of
$\lambda^+$ such that for all $\delta\in C$ limit, there is $a\in \PP_\delta$
such that
\[
\sup(a\cap X)=\delta.
\]
\end{description}
\end{Definition}

\begin{Fact}[D\v zamonja-Shelah]\label{545}\cite{DjSh 545}
If $\aleph_0<\kappa=\cf(\lambda)<\lambda$, \underline{then}
$\clubsuit^\ast_{-\lambda}(\lambda^+)$.
\end{Fact}

See the discussion below for a related
result of Shelah from \cite{Sh 667}.

Our Theorem \ref{general} can be understood as an extension of the
theorem from \cite{DjSh 545} to the successor of
a regular $\kappa$ for some $\kappa\ge\aleph_3$,
with the exception that our
guessing has less guesses at each $\delta$ (just one), but the guessing is obtained
stationarily often, as opposed to club often.
Also note
that the result in Theorem \ref{general}
is in some sense complementary to the ``club guessing"
results of Shelah (\cite{Sh 365} for example)
because here we are guessing unbounded
subsets of $\lambda^+$ which are not necessarily clubs, but on the other hand,
there are limitations on the cofinalities.

In \S\ref{druga}, we investigate successors of singulars. On the
one hand, we can hope to improve or at least
modify in a non-trivial way the above result from
\cite{DjSh 545}. If $\lambda$ is a strong limit singular,
\underline{and} $2^\lambda=\lambda^+$, it has already been
done by Shelah in \cite{Sh 667}:

\begin{Fact}[Shelah]\label{667}\cite{Sh 667} Suppose that $\lambda$ is a   
strong limit singular with $2^\lambda=\lambda^+$, and $S$ is a stationary
subset of $S^{\lambda^+}_{\cf(\lambda)}$.

\underline{Then} there is a sequence
\[
\left\langle\langle 
\bar{\alpha}^\delta=\alpha_{\delta,i}:\,i<\cf(\lambda)\rangle
:\,\delta\in S^{\lambda^+}_{\cf(\lambda)}\right\rangle
\]
such that $\bar{\alpha}^\delta$ increases to $\delta$, and for every
$\theta<\lambda$ and $f\in {}^{\lambda^+}\theta$, there
are stationarily many $\delta$ such that
\[
(\forall^\ast i) \,[f(\alpha_{\delta,2i})=f(\alpha_{\delta,2i+1})].
\]
\end{Fact}

Here, the quantifier $\forall^\ast i$ means 
``for all but $<\cf(\lambda)$ many". For more on guessing of unbounded sets,
see \cite{Sh -e}.

In \S\ref{druga}, we show that to a large extent the assumption
that $2^\lambda=\lambda^+$, is necessary above. See Theorem \ref{negth}.

We finish this introduction by recalling some notation and facts which will be
used in the following sections.

\begin{Notation}\label{oznake}
\begin{description}
\item{(1)}
Suppose that $\kappa=\cf(\kappa)<\delta$. We let
\[
S^\delta_\kappa\deq\{\alpha<\delta:\,\cf(\alpha)=\kappa\}.
\]
\item{(2)} Suppose that $C\subseteq \alpha$. We let
\[
\acc(C)\deq\{\beta\in C:\,\beta=\sup(C\cap \beta)\},
\]
and $\nacc(C)\deq C\setminus\acc(C)$, while $\lim(C)\deq
\{\delta<\alpha:\,\delta=\sup(C\cap\delta)\}$.
\end{description}
\end{Notation}

\begin{Definition}\label{square} Suppose that $\lambda\ge\aleph_1$
and $\gamma$ is an ordinal, while $A\subseteq \lambda^+$.
For $S\subseteq\lambda^+$, we say that $S$ {\em has a square
of type} $\le\gamma$ {\em nonaccumulating in} $A$ iff there
is a sequence $\langle e_\alpha:\,\alpha\in S\rangle$ such that
\begin{description}
\item{(i)} $\beta\in e_\alpha\implies \beta\in S\,\,\&\,\,e_\beta=e_\alpha
\cap\beta$,
\item{(ii)} $e_\alpha$ is a closed subset of $\alpha$,
\item{(iii)} If $\alpha\in S\setminus A$, then $\alpha=\sup(e_\alpha)$
(so $\nacc(e_\alpha)\subseteq A$ for all $\alpha\in S$),
\item{(iv)} $\otp{e_\alpha}\le\gamma$.
\end{description}
\end{Definition}

\begin{Fact}[Shelah]\label{351} [\cite{Sh 351}\S4, \cite{Sh g}III\S2]
Suppose that 
\[
\lambda=\cf(\lambda)=\theta^+>\theta=\cf(\theta)>\kappa=\cf(\kappa).
\]
Further suppose that $S\subseteq S^\lambda_\kappa$ is
stationary. \underline{Then} there is $S_1\subseteq\lambda$ on which there
is a square of type $\le\kappa$, nonaccumulating
on $A$=the successor ordinals, and $S_1\cap S$ is stationary.
\end{Fact}

\begin{Remark}\label{zabiljezlje} In the proof of Fact \ref{351} we can replace
$A$=the successor ordinals with $A=S^\lambda_\sigma$ for
any $\sigma=\cf(\sigma)
<\kappa$.
\end{Remark}

\begin{Definition}[Shelah]{\cite{Sh g}}\label{gl}
Suppose that $\delta<\lambda$ and
$e\subseteq \delta$, while $E\subseteq \lambda$. We define
\[
\gl(e,E)\deq\{\sup(\alpha\cap E):\,\alpha\in e\,\,\&\,\,\alpha>\min(E)\}.
\]
\end{Definition}

\begin{Observation} Suppose that $e$ and $E$ are as in Definition \ref{gl},
and both $e$ and $E\cap\delta$ are clubs of $\delta$. \underline{Then},
we have that
$\gl(e,E)$ is a club of $\delta$ with $\otp{\gl(e,E)}\le\otp{e}$.

If $e$ is just closed in $\delta$, and $E\cap\delta$
is a club of $\delta$ \underline{then} $\gl(e,E)$
is closed and $\otp{\gl(e,E)}\le\otp{e}$.
\end{Observation}

\begin{Fact}[Shelah]\label{355} \cite{Sh 365} Suppose that
$\cf(\kappa)=\kappa<\kappa^+<\cf(\lambda)=\lambda$. Further suppose that
$S\subseteq S^{\lambda}_{\kappa}$ is stationary and
$\langle e_\delta:\,\delta\in S\rangle$ is a sequence
such that each $e_\delta$ is a club of $\delta$.
\underline{Then}
there is a club $E^\ast$ of $\lambda$ such that the sequence
\[
\bar{c}=\langle c_\delta\deq\gl(e_\delta,E^\ast):\,\delta\in S
\cap E^\ast\rangle
\]
has the property that for every club $E$ of $\lambda$, there are stationarily
many $\delta$ such that $c_\delta\subseteq E$.
\end{Fact}

\begin{Observation}\label{both} Suppose that $\cf(\kappa)=\kappa<\kappa^+
<\lambda$ and $\lambda$ is a successor cardinal. Further assume that
$S_1\subseteq S^\lambda_\kappa$ is stationary, while $A=S^\lambda_\sigma$
for some $\sigma=\cf(\sigma)<\kappa$, possibly $\sigma=1$.
\underline{Then} there is stationary $S_2\subseteq S_1$
and a square
$\langle e_\delta:\,\delta\in S_2\rangle$ of
type $\le\kappa$ nonaccumulating in $A$,
such that each $e_\delta$ is a set of limit ordinals and $S_1\cap S_2$ is
stationary, while
\[
E\mbox{ a club of }\lambda\implies\{\delta\in S_1\cap S_2:\,e_\delta
\subseteq E\}\mbox{ is stationary}.
\]
\end{Observation}

[Why? The proof of this can be found in [\cite{Sh g}, III, \S2], but as it
easily follows from the Facts we already quoted, we shall give a proof.
By Fact \ref{351} and Remark \ref{zabiljezlje}, there is $S_3\subseteq \lambda$
with a square $\langle e_\alpha:\,\alpha\in S_3\rangle$ of type $\le
\kappa$
nonaccumulating in $A$, and such that $S_1\cap S_3$ is stationary. By Fact
\ref{355}, there
is  club $E^\ast$ of $\lambda$ as in the conclusion of Fact \ref{355},
with $S_3\cap S_1$ in place of $S$.
Now, letting
\[
S_2\deq\{\sup(\alpha\cap E^\ast):\,
\alpha\in \bigcup_{\delta\in S_3}e_\delta\cup\{\delta\}\,\,\&\,\,
\alpha>\min(E^\ast)\}\cap S_3,
\]
and for $\delta\in S_2$, letting
$c_\delta\deq\gl(e_\delta, E^\ast)$,
we observe that $S_2\cap S_1$ is stationary 
(as 
$S_2\cap S_1\supseteq S_1\cap S_3\cap\acc(E^\ast)$),
and $\langle c_\delta:\,
\delta\in S_2\rangle$ is a square of type $\le\kappa$
nonaccumulating in $A$, while
\[
E\mbox{ a club of }\lambda\implies\{\delta\in S_1\cap S_2:\,c_\delta
\subseteq E\}\mbox{ is stationary}.]
\]

\begin{Notation} Reg stands for the class of regular cardinals.
\end{Notation}

\section{A $ZFC$ version of $\clubsuit$}\label{prva}

\begin{Theorem}\label{general}
(1) Suppose that 
\begin{description}
\item{(a)}
$\lambda=\cf(\lambda)>\kappa=\cf(\kappa)>
\theta=\cf(\theta)\ge \aleph_1$.
\item{(b)} $S^\ast\subseteq S^{\lambda^+}_\theta$ is stationary,
moreover 
\[
S_1\deq\{\delta<\lambda^+:\,\cf(\delta)=\kappa
\,\,\&\,\,S^\ast\cap\delta\mbox{ is stationary}\}
\]
is stationary.
( e.g. $S^\ast=S^{\lambda^+}_{\theta}$.)
\end{description}
\underline{Then} there is a stationary $S'\subseteq S^\ast$ and
$\langle E_\delta:\,\delta\in S'\rangle$ such that
\begin{description}
\item{(i)} $E_\delta$ is a club of $\delta$ with $\otp{E_\delta}<
\lambda^\omega\cdot\kappa$,
\item{(ii)} for every unbounded $A\subseteq S^{\lambda^+}_\lambda$,
for stationarily many $\delta\in S'$, we have
\[
\delta=\sup(A\cap \nacc(E_\delta)).
\]
\end{description}

{\noindent (2)} If above we allow $\theta=\aleph_0$, but request
$\lambda\ge 2^{\aleph_0}$, the conclusion of (1) remains true.
\end{Theorem}

\begin{Remark} We explain why
the assumption that $S_1$ is stationary in item (b) above
implies that $S^\ast$ is stationary.
Notice that $\lambda^+>\kappa^+$, hence club guessing holds between
$\lambda^+$ and $\kappa$. In fact, by Fact \ref{355}, we can
assume that this is exemplified by a sequence $\langle c_\delta:\,
\delta\in S\subseteq S_1\rangle$. Now suppose that
$C$ is a club of $\lambda^+$, and let $E\deq \acc(C)$. Let $\delta\in S$
be such that $c_\delta\subseteq E$. Hence, $c_\delta$ is a club of $\delta$,
so $c_\delta\cap S^\ast\neq\emptyset$, implying that $C\cap S^\ast\neq
 \emptyset$.
 
 Obviously, $S_1$ being stationary is a necessary condition for our conclusion,
 as if $S_1$ were to be non-stationary, we could assume $S'\subseteq S^\ast
 \setminus S_1$, and $E_\delta\cap S^\ast=\emptyset$ for $\delta\in S'$. This
 would be a contradiction with (ii) above when $A=S^\ast$.
\end{Remark}

\begin{Proof} Let $S_0\deq S^{\lambda^+}_\theta$ and let $A^\ast\deq 
S^{\lambda^+}_{\aleph_1}$.

By Observation \ref{both} with $\lambda^+$ in place of
$\lambda$ and $A^\ast$ in place of $A$, there is a $S_2\subseteq
S^{\lambda^+}_{\le\kappa}$ such that there is a square
$\bar{e}=\langle e_\delta:\,\delta\in S_2\rangle$ of type $\le \kappa$
nonaccumulating in $A^\ast$,
the set $S_1\cap S_2$ is stationary, and, moreover, for every $E$ a club
of $\lambda^+$, the set $\{\delta\in S_1\cap S_2:\,e_\delta\subseteq E\}$
is stationary. [Why can we assume that $S_2\subseteq
S^{\lambda^+}_{\le\kappa}$? As if $\alpha\in S_2$ and $\alpha=\sup(e_\alpha)$,
we have that $\cf(\alpha)\le\card{e_\alpha}\le\kappa$ and if $\alpha\in S_2$
and $\alpha>\sup(e_\alpha)$ we have that $\alpha\in A^\ast$, hence
$\cf(\alpha)=\aleph_1\le\kappa$.

Let 
$S'\deq S^\ast \cap S_2$, so be stationary. [Why? As otherwise there is a
club $C$ of $\lambda^+$ with $S'\cap C=\emptyset$. Let $E\deq\acc(C)$ and
let $\delta\in S_1\cap S_2\cap E$ be such that $e_\delta\subseteq E$.
As $\cf(\delta)=\kappa\neq\aleph_1$, we have that $\delta\notin A^\ast$,
and hence $e_\delta$ is a club of $\delta$. On the other hand, $S^\ast\cap
\delta$ is stationary in $\delta$, hence $e_\delta\cap S^\ast\neq\emptyset$,
a contradiction with $e_\delta\subseteq C$, as $e_\delta\subseteq S_2$
by the definition of a square (see Definition \ref{square}). So any point in
$e_\delta
\cap S^\ast$ is in $C\cap S'$, contrary to the choice of $C$.]

\begin{Claim}\label{1A} There is a function $g:\,S'\to\omega$ such that
for every club $E$ of $\lambda^+$, there are stationarily many
$\delta\in S_1\cap S_2$ such that $e_\delta\subseteq E$ and
\[
(\forall n<\omega)[E\cap\delta\cap g^{-1}(\{n\})
\mbox{ is stationary in }\delta].
\]
\end{Claim}

\begin{Proof of the Claim} 
For $\delta\in S'$, we choose a sequence
$\bar{\xi}^\delta
=\langle \xi_{\delta,i}:\,i<\theta\rangle$ increasing with limit $\delta$,
and such that $\xi_{\delta,i}\in e_\delta$ and $\otp{e_{\xi_{\delta,i}}}$
depends only on $i$ and $\otp{e_\delta}$, but not on $\delta$.
[The point is of course that $\otp{e_\delta}$ is in general larger than
$\theta$.]
For each $i<\theta$, we define a function
$h_i:\,S'\into \kappa$ by letting 
\[
h_i(\delta)\deq \otp{e_{\xi_{\delta,i}}}.
\]

\begin{Subclaim}\label{A} For each $\delta\in S_1\cap S_2$ we can find
$i(\delta)<\theta$ such that with $i=i(\delta)$,
\[
A^\delta_i\deq\left\{\beta\in e_\delta:\,\{\gamma\in e_\delta\cap S':\,
\xi_{\gamma,i}=\beta\}\mbox{ is stationary }\right\}
\]
is unbounded in $\delta$.
\end{Subclaim}

\begin{Proof of the Subclaim} If this fails
for some $\delta\in S_1\cap S_2$, then each $A^\delta_i$
for $i<\theta$ is bounded in $\delta$. As $\theta<\kappa=\cf(\delta)$,
we have $\beta^\ast\deq\sup_{i<\theta}\sup (A^\delta_i)<\delta$.
As $\delta\in S_1$, we have that $e_\delta\cap S^\ast$ is stationary in $
\delta$. For every $\gamma\in e_\delta\cap S^\ast\setminus \beta^\ast$, we in
particular have that $\gamma\in S'$, so $\bar{\xi}^\gamma$
is defined. Hence, for such $\gamma$ there is $i_\gamma< \theta$
such that $\gamma>\xi_{\gamma,i_\gamma}>\beta^\ast$. By Fodor's
Lemma, there is $\xi^\ast$ such that for stationarily many $\gamma$
we have $\xi_{\gamma,i_\gamma}=\xi^\ast$, and applying the same lemma
again, we can without loss of generality assume that for some $i^\ast<\theta$
we have $i_\gamma=i^\ast$ for stationarily many $\gamma$ for which
$\xi_{\gamma,i_\gamma}=\xi^\ast$. But then $\xi^\ast>\beta^\ast$
and yet $\xi^\ast\in A^\delta_{i^\ast}$, a contradiction.
$\eop_{\ref{A}}$
\end{Proof of the Subclaim}

\begin{Subclaim}\label{B} For some $i(\ast)<\theta$, the set
\[
S^{\ast\ast}\deq\{\delta\in S_1\cap S_2:\,i(\delta)=i(\ast)\}
\]
is stationary, and $\bar{e}\rest S^{\ast\ast}$ still guesses clubs
of $\lambda^+$.
\end{Subclaim}

\begin{Proof of the Subclaim} Otherwise, for each $i<\theta$ such that
\[
T_i\deq\{\delta\in S_1\cap S_2:\,i(\delta)=i\}
\]
is stationary, there is a club $C_i$ of $\lambda^+$ such that for no $\delta
\in T_i$ do we have $e_\delta\subseteq C_i$. Let $C\deq\bigcap\{C_i:\,
T_i\mbox{ stationary }\}$,
hence a club of $\lambda^+$, and let $E$ be a club of $\lambda^+$
such that
\[
[i<\theta\,\,\&\,\,T_i\mbox{ not stationary}]\implies T_i\cap E=\emptyset.
\]
Let $\delta\in S_1\cap S_2$ be such that $e_\delta\subseteq\acc(E\cap C)$.
Hence $\delta\in E\cap T_{i(\delta)}$, so $T_{i(\delta)}$
is stationary. On the other hand, we have $e_\delta\subseteq C_{i(\delta)}$,
a contradiction.
$\eop_{\ref{B}}$
\end{Proof of the Subclaim}

Now notice that $\kappa>\aleph_1$, so club guessing holds between
$\kappa$ and $\aleph_0$, i.e. there is a sequence
$\langle w_\zeta:\,\zeta\in S^\kappa_{\aleph_0}\rangle$ such that
$w_\zeta\subseteq\zeta$ and $\otp{w_\zeta}=\omega$ for each $\zeta$, while
for every club $C$ of $\kappa$, there are stationarily many $\zeta$
with $w_\zeta\subseteq C$. Let $W\deq\{w_\zeta:\,\zeta\in S^\kappa_{\aleph_0}
\}$. For $\beta\in S'$, let
\[
w^\zeta_\beta\deq\{\gamma\in e_\beta:\,\otp{e_\beta\cap\gamma}\in w_\zeta\}.
\]
For each $\zeta\in S^\kappa_{\aleph_0}$, we define $g_\zeta:\,S'
\into\omega$ by letting
\[
g_\zeta(\gamma)\deq\otp{h_{i(\ast)}(\gamma)\cap w_\zeta }.
\]
If some $g_\zeta$ is as required in Claim \ref{1A}, then we are done.
Otherwise, for each $\zeta$ there is a club $E_\zeta$ of $\lambda^+$
with
\[
[\delta\in S^{\ast\ast}\cap E_\zeta\,\&\,e_\delta\subseteq E_\zeta]
\implies\,(\exists n\deq n_{\delta,\zeta})\,[E_\zeta\cap\delta\cap
g^{-1}_\zeta(\{n\})\mbox{ is non-stationary in }\delta].
\]
Let $E\deq\bigcap_{\zeta<\kappa}E_\zeta$. 
Let $\delta^\ast\in S^{\ast\ast}\cap E$ be such that $e_{\delta^\ast}
\subseteq E$. As $\delta^\ast\in S^{\ast\ast}$, we have that
$A^{\delta^\ast}_{i(\ast)}$ is unbounded in $\delta^\ast$. Let
$C_{\delta^\ast}\deq\lim(A^{\delta^\ast}_{i(\ast)})\cap e_{\delta^\ast}$,
hence a club of $\delta^\ast$. As $\cf(\delta^\ast)=\kappa$, and so
$\otp{e_{\delta^\ast}}=\kappa$, we have that
\[
C^-_{\delta^\ast}\deq\{\otp{\gamma\cap e_{\delta^\ast}}:\,\gamma\in
C_{\delta^\ast}\}
\]
is a club of $\kappa$. Hence there is $\zeta^\ast\in S^\kappa_{\aleph_0}$
such that $w_{\zeta^\ast}\subseteq C^-_{\delta^\ast}$. Let
$\{\xi_0,\xi_1,\xi_2,\ldots\}$ be the increasing enumeration of
$w_{\zeta^\ast}$, and stipulate $\xi_{-1}=0$. For $l<\omega$, let
$\gamma_l$ be the unique $\gamma\in C_{\delta^\ast}$ such that
$\otp{\gamma_l\cap e_{\delta^\ast}}=\xi_l$, and let $\gamma_{-1}
\deq\min(C_{\delta^\ast})$. Hence for every $l<\omega$ we have
\[
A^{\delta^\ast}_{i(\ast)}\cap [\gamma_{l-1},\gamma_l)\neq\emptyset.
\]
Hence for some $\beta_l\in [\gamma_{l-1},\gamma_l)$, we have that
\[
\{\gamma\in e_{\delta^\ast}\cap S':\,\xi_{\gamma,i(\ast)}=\beta_l\}
\]
is stationary in $\delta^\ast$. This means that for each $l<\omega$,
the set $g^{-1}_{\zeta^\ast}(\{l\})\cap e_{\delta^\ast}\cap S'$
is stationary in $\delta^\ast$, a contradiction.
$\eop_{\ref{1A}}$
\end{Proof of the Claim}

{\em Continuation of the Proof of Theorem \ref{general}}.
We now choose 
$\bar{c}=\langle c_\alpha:\,\alpha <\lambda^+\rangle$, so that
\begin{description}
\item{($\alpha$)} For every $\alpha$, we have that
$c_\alpha$ is a club of $\alpha$ with $\otp{c_\alpha}\le\lambda$,
and for $\alpha$ a limit
ordinal $\beta\in \acc(c_\alpha)\implies\cf(\beta)<\lambda$, while
if $\alpha=\beta+1$, then $c_\alpha=\{\beta\}$.

\item{($\beta$)} If $\delta\in S_2$, then $c_\delta\supseteq e_\delta$,

\item{($\gamma)$} If $\delta\in S_2$ and $\sup(e_\delta)=\delta$, then
$c_\delta=e_\delta$,

\end{description}



Now for any limit $\delta<\lambda^+$ we choose by induction on $n<\omega$
a club $C^n_\delta$ of $\delta$ of order type $\le\lambda^{n+1}$,
using the following algorithm:

Let $C^0_\delta\deq c_\delta$.
Let
\[
C^{n+1}_\delta\deq C^n_\delta \cup\{\alpha:\,(\exists\beta\in
\nacc(C^n_\delta))\,[\sup(\beta\cap C^n_\delta)<\alpha
<\beta \,\,\&\,\,\alpha\in c_\beta]\}.
\]

\begin{Note}\label{Cn}
\begin{description}
\item{(1)} The above algorithm really gives $C^n_\delta$
which is a club of $\delta$
with
\[
\otp{C^n_\delta}\le\lambda^{n+1}.
\]
If $\delta\in S_2$ and $\sup(e_\delta)=\delta$, then
$\otp{C^\delta_n}\le\lambda^n\cdot\kappa$.

[Why? We prove this by induction on $n$. It is clearly true for $n=0$.
Assume its truth for $n$. Clearly $C^{n+1}_\delta$ is unbounded in $\delta$,
let us show that it is closed. Suppose $\alpha=\sup(C^{n+1}_\delta\cap
\alpha)<\delta$. If $\alpha=\sup(C^{n}_\delta\cap
\alpha)$, then $\alpha\in C^n_\delta\subseteq C^{n+1}_\delta$ by the induction
hypothesis. So, assume
\[
\alpha^\ast\deq\sup(C^{n}_\delta\cap
\alpha)<\alpha
\]
and
$\alpha\notin C^n_\delta$. Let $\langle\alpha_i:\,i<\cf(\alpha)\rangle$ 
be an increasing to $\alpha$
sequence in $(\alpha^\ast,\alpha)\cap C^{n+1}_\delta$. Hence for every $i$
there is $\beta_i\in \nacc(C^n_\delta)$ such that $\alpha_i\in c_{\beta_i}$
and $\sup(C^n_\delta\cap\beta_i)<\alpha_i$. As $\sup(C^n_\delta\cap\alpha)=
\alpha^\ast<\alpha_i$ and $\alpha\notin C^n_\delta$, we have $\beta_i>\alpha$,
for every $i$,
as $\beta_i\in\nacc(C^n_\delta)$. 
Suppose that $i\neq j$ and $\beta_i<\beta_j$. Hence
$\sup(C^n_\delta\cap\beta_j)\ge\beta_i>\alpha_j$, a contradiction.
So, there is $\beta$ such that $\beta_i=\beta$ for all $i$, hence
$\{\alpha_i:\,i<\cf(\alpha)\}\subseteq c_\beta$. As $c_\beta$ is closed,
and $\alpha<\beta$,
we have $\alpha\in c_\beta$, and by the definition of $C^{n+1}_\delta$
we have $\alpha\in C^{n+1}_\delta$.

As for every $\beta$ we have $\otp{c_\beta}\le\lambda$, and by the induction
hypothesis $\otp{C^n_\delta}\le\lambda^{n+1}$, we have $\otp{C^{n+1}_\delta}
\le\lambda^{n+2}$.

Similarly, if $\delta\in S_2$ and $\sup(e_\delta)=\delta$, clearly
$\otp{C^n_\delta}\le\lambda^n\cdot \kappa$.]

\item{(2)} For every $n$, we have $\acc(C^n_\delta)\setminus
\bigcup_{m<n}C^m_\delta\subseteq
S^{\lambda^+}_{<\lambda}$.

[Why? Again by induction on $n$. For $n=0$ it follows as $\otp{c_\delta}
\le\lambda$. Suppose this is true for $C^n_\delta$. The analysis from the proof
of
(1) shows that for $\alpha\in \acc(C^{n+1}_\delta)\setminus
C^n_\delta$, there is $\beta$ such that $\alpha\in \acc(c_\beta)$, hence
$\cf(\alpha)<\lambda$.]

\item{(3)} For every limit $\delta<\lambda^+$, we have $S^\delta_\lambda=
\bigcup_{n<\omega}\nacc(C^n_\delta)\cap S^\delta_\lambda$.

[Why? Fix such $\delta$ and let $\alpha\in S^\delta_\lambda$.
By item (2), it suffices to show that $\alpha\in C^n_\delta$ for some $n$.
Suppose not, so let $\gamma_n\deq\min(C^n_\delta\setminus\alpha)$ for $n<\omega$.
Hence $\langle\gamma_n:\,n<\omega\rangle$ is a non-increasing
sequence of ordinals $>\alpha$, and so there is $n^\ast$ such that
$n\ge n^\ast\implies\gamma_n=\gamma_{n^\ast}$. In particular we have that
$\gamma_{n^\ast}\in \nacc(C^{n^\ast}_\delta)$. Let $\beta\in c_{\gamma_{n^\ast}}
\setminus\alpha$. Hence $\sup(\gamma_{n^\ast}
\cap C^{n^\ast}_\delta)<\alpha\le\beta
<\gamma_{n^\ast}$. By the definition of $C^{n^\ast+1}_\delta$, we have
$\beta\in C^{n^\ast+1}_\delta$, a contradiction.]
\end{description}
\end{Note}

Now for each $\delta\in S'$ we define
\[
E_\delta\deq e_\delta\cup\bigcup\{C^{g(\delta)}_\alpha
\setminus\sup(e_\delta\cap\alpha):\,
\alpha\in\nacc(e_\delta)\}.
\]
Note first that $E_\delta$ is a club of $\delta$, for $\delta\in S'$.

[Why? Clearly, $E_\delta$ is unbounded. Suppose $\gamma=\sup(E_\delta
\cap\gamma)<\delta$.
Without loss of generality we can assume $\gamma\notin e_\delta$.
Let $\gamma^\ast\deq \sup(e_\delta\cap\gamma)<\gamma$.
For every $\beta\in E_\delta\cap (\gamma^\ast,\gamma)$, there is
$\alpha_\beta\in\nacc(e_\delta)\cap S'$ such that
$\beta\in C^{g(\delta)}_{\alpha_\beta}\setminus
\sup(e_\delta\cap\alpha_\beta)$.
By the choice of $\gamma^\ast$, every such $\alpha_\beta>\gamma$.
Suppose that $\beta_1\neq\beta_2\in E_\delta\cap (\gamma^\ast,\gamma)$
and $\alpha_{\beta_1}<\alpha_{\beta_2}$. Hence $\sup(e_\delta\cap
\alpha_{\beta_2})\ge\alpha_{\beta_1}$, a contradiction. So all
$\alpha_\beta$ are a fixed $\alpha$. Hence $\gamma<\alpha$ is a limit point of
$C^{g(\delta)}_\alpha$, and we are done, as $C^{g(\delta)}_\alpha$ is closed.]

Also note that $\otp{E_\delta}<\lambda^\omega\cdot \kappa$.

Suppose that $A\subseteq S^{\lambda^+}_\lambda$ is unbounded
and it exemplifies that $\langle E_\delta:\,\delta\in S'\rangle$
fails to satisfy the requirements of Theorem
\ref{general}. Hence there is a club  $E$ of
$\lambda^+$ such that
\[
\delta\in E\cap S'\implies \sup(A\cap\nacc(E_\delta))<\delta.
\]
Let $E^\ast\deq\acc(E)\cap\{\delta:\,\delta=\sup(A\cap\delta)\}$, hence a
club of $\lambda^+$. Let $\delta^\ast\in S_1\cap S_2\cap E^\ast$ be such that
$e_{\delta^\ast}\subseteq E^\ast$ and for all $n<\omega$,
the set $\delta^\ast\cap g^{-1}(\{n\})$ is stationary in $\delta^\ast$.

For $\alpha\in \nacc(e_{\delta^\ast})$ we have that $A\cap \alpha$
is unbounded in $\alpha$.
Now we use Note \ref{Cn}(3).
As $A\subseteq S^{\lambda^+}_\lambda$ we have $A\cap\alpha
\subseteq S^\alpha_\lambda$. So $A\cap\alpha=\bigcup_{n<\omega}
\nacc(C^n_\alpha)\cap
A\cap\alpha$, by the above mentioned Note. As $\alpha\in
\nacc(e_{\delta^\ast})$
and $\nacc(e_{\delta^\ast}) \subseteq S^{\lambda^+}_{\aleph_1}$, there is
$n<\omega$ such that $A\cap\nacc(C^n_\alpha)$ is unbounded in $\alpha$.
Let $n^\ast(\alpha)$ be the smallest such $n$. There is $n^\ast$
such that
\[
\sup\{\alpha\in\nacc(e_{\delta^\ast}):\,n^\ast(\alpha)=n^\ast\}=\delta^\ast,
\]
as $\cf(\delta^\ast)>\aleph_0$. Let
\[
e\deq\left\{\beta\in\acc(e_{\delta^\ast}):\,\beta=\sup\{\alpha\in\beta\cap
\nacc(e_{\delta^\ast}):\,n^\ast(\alpha)=n^\ast\}\right\},
\]
hence $e$ is a club of $\delta^\ast$. By the choice of $\delta^\ast$,
the set $g^{-1}(\{n^\ast\})
\cap\delta^\ast$ is stationary in $\delta^\ast$. So, there is
$\beta\in e$ such that $g(\beta)=n^\ast$. In particular,
$\beta\in S'$.
For every $\alpha\in \nacc(e_\beta)$ such that
$n^\ast(\alpha)=n^\ast$, we have
that $A\cap \nacc(C^{n^\ast}_\alpha)$ is unbounded
in $\alpha$. However,
\[
C^{n^\ast}_\alpha\setminus\sup(\alpha\cap e_{\beta})=
E_{\beta}\cap [\sup(\alpha\cap e_\beta),\beta),
\]
hence $\nacc(C^{n^\ast}_\alpha)\setminus\sup(\alpha\cap e_\beta)
\subseteq\nacc(E_\beta)$.
Now, on the one hand $\beta\in E\cap S'$, so $\alpha^\ast\deq\sup(A\cap
\nacc(E_\beta))<\beta$, but on the other hand, the set of
$\alpha\in\nacc(e_\beta)$ with $n^\ast(\alpha)=n^\ast$
is unbounded in $\beta$,
hence there is $\gamma\in A
$ with $\gamma\in E_{\beta}\setminus
\alpha^\ast$. As $\gamma\in A$, we have $\cf(\gamma)=\lambda$,
hence $\gamma\in \nacc(E_{\beta})$, a contradiction.

{\noindent (2)} The statement of Claim \ref{1A}
is true even when $\kappa=\aleph_1$, but if we assume that $\lambda
\ge 2^{\aleph_0}$. Namely, under these assumptions, we have that
$\theta=\aleph_0$, so there is $W\subseteq\{w\subseteq\omega_1:\,\otp{w}
=\omega\}$ such that $\card{W}\le\lambda$, and for every club $C$
of $\omega_1$, for some $w\in W$, we have $w\subseteq C$. Now we can
just repeat the proof of Claim \ref{1A}, using $W$ we have just defined.
$\eop_{\ref{general}}$
\end{Proof}

\section{A negation of guessing}\label{druga}

\begin{Theorem}\label{negth}
Assume that there is a supercompact cardinal.
\underline{Then}

\begin{description}
\item{(1)} It is consistent
that there is $\lambda$ a strong limit singular
of cofinality $\aleph_0$, such that $2^\lambda>\lambda^+$ and
\begin{description}
\item{$(\ast)$} There is a function $f:\,\lambda^+\into
\omega$ such that for every $\PP\subseteq [\lambda^+]^{\aleph_0}$ of
cardinality $<2^\lambda$, for some $X\in [\lambda^+]^{\lambda^+}$ we have
\begin{description}
\item{(i)} $(\forall \zeta<\omega)[\card{X\cap f^{-1}(\{\zeta\})}=\lambda^+]$,
\item{(ii)} If $a\in \PP$, then $\sup({\rm Rang}(f\rest(a\cap X))<\omega$.
\end{description}
\end{description}
\item{(2)} Moreover, in (1) we can replace $\aleph_0$ by any regular
$\kappa<\lambda$.
\end{description}
\end{Theorem}

\begin{Remark} So the theorem basically
states that no $\PP$ as above provides a 
guessing.
\end{Remark}

\begin{Proof} (1) We start with a universe in which
$\lambda$ is a supercompact cardinal and $2^\lambda=\lambda^+$ holds.
We extend the universe by Laver's
forcing
(\cite{Laver}), which makes the supercompactness of
$\lambda$ indestructible by any extension by a
$(<\lambda)$-directed-closed forcing. This forcing will preserve the fact that
$2^{\lambda}=\lambda^+$.
Let us call the so obtained universe $V$.

Now choose $\mu$ such that $\mu=\mu^\lambda>\lambda^+$.
By \cite{Baumgartner}, there is a $(<\lambda)$-directed-closed
$\lambda^{++}$-cc forcing notion
$P$, not collapsing $\lambda^+$ of size $\mu$ adding $\mu$
unbounded subsets $A_\alpha\,(\alpha<\mu)$ to $\lambda^+$ such that
\begin{description}
\item{$(\ast\ast)$}
$\alpha\neq \beta<\mu\implies\card{A_\alpha\cap A_\beta}<\lambda$.
\end{description}
In particular, in $V^P$ we have
$\lambda^+<2^{\lambda}= \mu$
(\cite{Baumgartner}, 6.1.), while
$\lambda$ is supercompact.
In $V^P$, let $Q$ be Prikry's forcing which does not collapse cardinals
and makes $\lambda$ singular with $\cf(\lambda)=\aleph_0$,
\cite{Prikry}. As this forcing does not add bounded
subsets to $\lambda$, in the extension $\lambda$ is a strong limit singular
and clearly
satisfies $2^\lambda=\mu$. In
$V^{P\ast\name{Q}}$ we have $(\ast\ast)$. We now work in $V^{P\ast\name{Q}}$.

Let $\lambda=\sum_{\zeta<\omega}\lambda_\zeta$ where each
$\lambda_\zeta<\lambda$ is
regular. Let
$\chi$ be large enough regular and
$M\elementary (\HH(\chi),\in)$ with $\card{\card{M}}=\lambda^+$
such that $\lambda^+\subseteq M$ and
$\langle A_\alpha:\,\alpha<\mu\rangle,
\langle \lambda_\zeta:\,\zeta<\omega\rangle\in M$. We list
$\bigcup_{\zeta<\omega}([\lambda^+]^{\lambda_\zeta}\cap M)$ as
$\{b_i:\,i<\lambda^+\}$.

We define $f:\,\lambda^+\into\omega$ by
$f(i)=\zeta$ iff $\card{b_i}=\lambda_\zeta$. For $\alpha<\mu$, let
$X_\alpha\deq\{i:\,b_i\subseteq A_\alpha\}$.

Now suppose that $\PP\subseteq [\lambda^+]^{\aleph_0}$ is of cardinality
$<2^\lambda\le\mu$, we shall look for $X$ as required in $(\ast)$.

If $\alpha<\mu$ is such that $X_\alpha$ fails to serve as $X$, then
at least one of the following two cases must hold:

\underline{Case 1}. For some $\zeta<\omega$ we have
$\card{\{i:\,b_i\subseteq
A_\alpha\,\,\&\,\,\card{b_i}=\lambda_\zeta\}}<\lambda^+$,
or

\underline{Case 2}. For some $a\in\PP$ we have
$\sup({\rm Rang}(f\rest (a\cap X_\alpha)))=\omega$.

Considering the second case, we shall show that for any $a\in \PP$,
there are $<\lambda$ ordinals $\alpha$ such that the second case
holds for $X_\alpha, a$. Fix an $a\in \PP$. If $\alpha<\mu$ is such
that Case 2 holds for $X_\alpha,a$, then
\[
\sup(\{\zeta:\,(\exists i\in a)[b_i\subseteq A_\alpha\,\,\&\,\,
\card{b_i}=\lambda_\zeta]\}=\omega.
\]
For $\zeta<\omega$
and $\alpha<\mu$ let $B_{\zeta}^\alpha\deq\{i\in a:\,b_i\subseteq A_\alpha
\,\,\&\,\,
\card{b_i}=\lambda_\zeta\}$. Notice that if $\alpha\neq\beta<\mu$ we have that
for some $\zeta_{\alpha,\beta}$ the intersection $A_\alpha\cap A_\beta$
has size $<\lambda_{\zeta_{\alpha,\beta}}$, hence for all $\zeta\ge
\zeta_{\alpha,\beta}$ we have $B_{\zeta}^\alpha\cap B_{\zeta}^\beta=
\emptyset$.

Let $A\deq\{\alpha:\,\mbox{ Case 2 holds for }a,\alpha\}$.
For every $\alpha\in A$, let 
\[
\bar{s_\alpha}\deq\langle \,B^\alpha_\zeta:\zeta<\omega\rangle,
\]
hence $\alpha\neq\beta\implies \bar{s_\alpha}\neq \bar{s_\beta}$.
Hence $\card{A}\le 2^{\aleph_0}<\lambda$.

Now note that if $\alpha<\mu$, then $A_\alpha\in [\lambda^+]^{\lambda^+}$.
For $\gamma_0<\lambda^+$ for some $\gamma_1<\lambda^+$ we have
$\card{A_\alpha\cap\gamma_1\setminus\gamma_0}=\lambda$. In $M$
we have a sequence $\langle c_\xi:\,\xi<\omega\rangle$
such that $\cup_\xi c_\xi=\gamma_1\setminus\gamma_0$ and $\card{c_\xi}=
\lambda_\xi$. For every $\zeta, \varepsilon<\omega$, for some
large enough $\xi<\omega$ we
have $\card{A_\alpha\cap c_\xi}\ge\lambda_\varepsilon$. But $[A_\alpha\cap
c_\xi]^{\lambda_\zeta}\subseteq \PP(c_\xi)\subseteq M$ (as
$\lambda^{<\lambda}<\lambda$), so
for some $i$ we have $b_i\subseteq A_\alpha\cap [\gamma_0,\gamma_1)$
and $\card{b_i}=\lambda_\zeta$. Hence
\[
\card{\{i:\,b_i\subseteq A_\alpha\,\,\&\,\,
\card{b_i}=
\lambda_\zeta\}}=\lambda^+,
\]
so Case 1 does not happen
for this (any) $\alpha$.

As we can find $\alpha <\mu$ such that Case 2 does not happen, we
are finished.

(2) Use Magidor's forcing from \cite{Magidor} in place of Prikry's
forcing in (1).
$\eop_{\ref{negth}}$
\end{Proof}

\eject

\vskip 5truecm

The address of the author for correspondence:

Mirna D\v zamonja

School of Mathematics

University of East Anglia

Norwich, NR4 7TJ

UK

Telephone: +44-1603-592981

Fax:+44-1603-593868

e-mail: M.Dzamonja@uea.ac.uk
\end{document}